\documentclass[aps,prd,twocolumn,superscriptaddress,nofootinbib,floatfix,amsmath,amssymb]{revtex4-2}

\usepackage{amsmath,amssymb,mathtools,bm}
\usepackage{graphicx}
\usepackage{booktabs}
\usepackage{array}
\usepackage{enumitem}
\usepackage{microtype}
\usepackage{xcolor}
\usepackage[colorlinks=true,linkcolor=blue,citecolor=blue,urlcolor=blue]{hyperref}

\newcommand{\R}{\mathbb{R}}
\newcommand{\sign}{\operatorname{sign}}
\newcommand{\tr}{\operatorname{tr}}

\newcommand{\argmin}{\operatorname*{arg\,min}}

\newcommand{\GD}{\mathrm{GD}}
\newcommand{\SGN}{\mathrm{SIGN}}

\begin{document}

\title{Adaptive Gradient Geometry with Memory:\\
Local Optimality, Trajectory Value, and the Limits of Greedy Selection}

\author{Lorenzo Pepa}
\email{pepalorenzo.work@gmail.com}
\affiliation{Hewlett Packard Labs, HPE Belgium, B-1831 Diegem, Belgium}
\author{Mohammad Ali Marashli}
\affiliation{Lule\aa{} University of Technology, Lule\aa{}, Sweden}
\author{Tommaso Dorigo}
\affiliation{Lule\aa{} University of Technology, Lule\aa{}, Sweden}
\affiliation{Istituto Nazionale di Fisica Nucleare -- Sezione di Padova, via F. Marzolo 8, 35131 Padova, Italy}
\author{Fredrik Sandin}
\affiliation{Lule\aa{} University of Technology, Lule\aa{}, Sweden}

\date{\today}

\begin{abstract}
In minimization problems, the direction in parameter space that lowers a cost function most in a single step need not be the best choice over several steps, because the choice also changes what the algorithm remembers. We investigate this effect for gradient descent, sign descent (which keeps only the signs of gradient components), and intermediate power-gradient directions, all combined with a correction that uses information about local curvature and the preceding step.

Examination of an example involving an exact three-dimensional positive-definite quadratic surface shows that the direction giving the larger immediate decrease in cost can lead to a higher cost after two steps. In numerical experiments on three families of smooth convex functions, an idealized two-step search shows that changing the gradient direction can be useful away from the nearly quadratic regime. When intermediate ways of weighting gradient components are allowed, rather than just gradient descent and sign descent, an intermediate first direction is preferred in about one third of the sampled states.

A diagnostic based on curvature changes along the
preceding step predicts the average preferred exponent well, even when the family being tested is excluded from the fit. But using that prediction to choose directions during minimization does not reliably outperform gradient descent with the same memory. We find that knowing which exponent is preferred on average does not tell us which choice will help a particular trajectory.

\end{abstract}
\maketitle

\section{Motivation and scope}

Many problems in experimental and computational physics require the minimization of a cost function which depends on several coupled parameters. Examples include fitting models to experimental measurements, calibrating detector responses, tuning reconstruction algorithms, and optimizing experimental configurations using fast surrogate models of expensive simulations. 

In each case, a gradient-based method uses the derivatives with respect to the individual parameters to choose a direction in parameter space: the numerical value of each of them depends on the units in which the corresponding parameter is expressed. Since physical parameters may have different units and characteristic scales, gradient components may have numerical magnitudes that differ by several orders of magnitude. In ordinary gradient descent, the step length may be limited by the coordinates with the largest gradient components, leaving little progress along other directions that could be equally important for reaching the minimum. A standard remedy is to express each parameter in units of a characteristic physical scale, so that a unit change in each rescaled coordinate then corresponds to a change of one characteristic scale in the underlying physical parameter. More generally, diagonal preconditioning assigns different weights to the parameter updates, using information about their typical sensitivities or curvatures. Such a choice can considerably improve convergence, although a fixed diagonal rescaling cannot generally eliminate correlations between parameters or accommodate curvature that changes along the minimization trajectory.

Even after suitable rescaling, the relative magnitudes of gradient components may vary considerably during minimization. How much should the optimizer rely on those magnitudes when choosing its next direction? Ordinary gradient descent uses them directly; sign descent retains only their signs. With intermediate power-gradient directions we can vary the relative weight of the components continuously.

If local curvature information is available, it seems sensible to choose the direction that is expected to give the largest decrease in cost. Many successful minimization algorithms, however, use information from earlier steps to avoid repeatedly moving along unfavorable directions. In this case, the direction chosen for one step also enters into the calculation of the next one; as a result, the best immediate step may turn out to be a poor choice in the long run.

In this work we study how changing the way gradient components are weighted interacts with the information retained from previous steps. In particular, we want to understand whether a direction that looks favorable at the present step remains a good choice once its consequences for the subsequent trajectory are taken into account. We also investigate whether useful changes of direction can be identified from curvature information available during minimization.

This document is organized as follows. Section~II introduces the gradient geometries considered in this study, and Sec.~III describes the curvature-memory correction. Section~IV presents an exact quadratic counterexample to locally greedy direction selection. Sections~V and VI
examine the value of adapting gradient geometry on
nonlinear cost functions and its relationship with changes in local curvature. Then in Section~VII we test whether that relationship can guide practical optimization rules. Related minimization methods are discussed in Sec.~VIII, followed by a discussion of our results in Sec.~IX and the conclusions in Sec.~X.

\section{Local geometry: gradient descent, SIGN, and power gradients}
\label{sec:local}

Let $f:\R^d\to\R$ be a smooth cost function, expressed in
a chosen set of dimensionless coordinates. We denote its
gradient and Hessian at the current point $x$ by
$g=\nabla f(x)$ and $H=\nabla^2 f(x)$, respectively.
For a small displacement $\delta$, the quadratic
approximation is
\begin{equation}
 f(x+\delta)\simeq f(x)+g^T\delta
                      +\frac12\delta^TH\delta.
 \label{eq:localquad}
\end{equation}

Consider a search direction $v$ satisfying $g^Tv<0$
and $v^THv>0$. Minimizing the quadratic approximation
along the ray $x+\alpha v$ gives
\begin{equation}
 \alpha_{\rm quad}=-\frac{g^Tv}{v^THv},
\end{equation}
and a predicted decrease
\begin{equation}
 D(v)=\frac{(g^Tv)^2}{2v^THv}.
 \label{eq:D}
\end{equation}
The numerator in Eq.~\eqref{eq:D} is the squared directional derivative for a displacement along $v$; the denominator contains the corresponding curvature. For a quadratic cost function the decrease is exact, while for a general smooth function it is only a local prediction.

The notion of steepest descent depends on how the
size of a displacement is measured. With the
Euclidean norm, the direction of greatest first-order
decrease is the negative gradient, $v=-g$. If instead
we constrain the largest absolute coordinate
displacement, using the $\ell_\infty$ norm, the
corresponding direction is $v=-\sigma$, where
$\sigma_i=\sign(g_i)$, with $\sign(0)=0$.
Thus gradient descent (GD) and sign descent (SIGN)
represent two different choices of descent geometry.

For GD, Eq.~\eqref{eq:D} becomes
\begin{equation}
 D_{\GD}=\frac{\|g\|_2^4}{2g^THg},
\end{equation}
whereas for SIGN,
\begin{equation}
 D_{\SGN}=\frac{\|g\|_1^2}{2\sigma^TH\sigma}.
\end{equation}
Their ratio is
\begin{equation}
 \boxed{
 R=\frac{D_{\SGN}}{D_{\GD}}
 =\frac{\|g\|_1^2\,g^THg}
        {\|g\|_2^4\,\sigma^TH\sigma}.
 }
 \label{eq:R}
\end{equation}
A value $R>1$ indicates that SIGN gives the larger
predicted decrease, while $R<1$ favors GD.

We can separate the factors controlling this ratio:
\begin{equation}
 R=G\,C\,B,
\end{equation}
where
\begin{equation}
 \begin{split}
 G&=\frac{\|g\|_1^2}{d\|g\|_2^2},\\
 C&=\frac{d\,g^THg}{\|g\|_2^2\tr H},\\
 B&=\frac{\tr H}{\sigma^TH\sigma}.
 \end{split}
\end{equation}
The first factor measures the distribution of gradient
magnitudes over the coordinates: $G$ approaches one
when their absolute values are similar. The second compares the curvature along the gradient with the mean curvature \(\operatorname{tr}H/d\), where \(\operatorname{tr}H\) denotes the trace of the Hessian, the sum of its diagonal elements.
The third reflects the
curvature encountered along the SIGN direction,
including the effects of off-diagonal Hessian terms.
Then we observe that for an isotropic Hessian $H=\lambda I$ and a gradient
with no vanishing components, $C=B=1$ and
$R=G\leq 1$. More generally, the same inequality
follows from the Cauchy--Schwarz inequality restricted
to the nonzero components of $g$. Thus SIGN cannot
improve the immediate quadratic-model decrease
over GD in an isotropic landscape. An advantage
for SIGN requires a suitable combination of gradient
structure and anisotropic curvature.

Between gradient descent and sign descent we can define intermediate prescriptions that retain some information about the relative sizes of the gradient components. We introduce a dimensionless exponent \(h\), ranging from zero to one, to control how strongly these magnitudes influence the search direction. The resulting family of directions is

\begin{equation}
 \boxed{
 q_i(h)=-\sign(g_i)|g_i|^h,\qquad 0\leq h\leq1.
 }
 \label{eq:power}
\end{equation}
For $h=0$ we define $q_i(0)=-\sign(g_i)$,
including when $g_i=0$, so that $q_i(0)=0$
for a vanishing gradient component.
Here $h=1$ gives GD, while $h=0$ gives SIGN.
Intermediate exponents preserve the ordering of
gradient magnitudes while reducing their relative
differences. Up to an overall normalization, the
direction $q(h)$ is a steepest-descent direction
associated with the $\ell_{1+1/h}$ norm for $h>0$,
with SIGN recovered in the $\ell_\infty$ limit.
This family is closely related to the Powerball
transformation~\citep{Yuan2019Powerball}.

The choice of coordinates matters here. Let $\theta_i$ denote a physical
parameter and write $\theta_i=c_i x_i$, where $c_i>0$
is a characteristic physical scale and $x_i$ is
dimensionless. If $g_i^{(\theta)}$ is the gradient
component with respect to $\theta_i$, applying GD in
the dimensionless coordinates produces physical
updates proportional to $-c_i^2g_i^{(\theta)}$.
SIGN instead gives physical updates proportional to
$-c_i\sign(g_i^{(\theta)})$. Their relative behavior
therefore depends on the parameter scales. Throughout our analysis, the descent geometry refers to the fixed, dimensionless coordinates in which the minimization is performed.

Equation~\eqref{eq:D} allows us to compare the candidate directions for a single step. If the algorithm retains information from earlier steps, however, our choice also affects the next search direction.

\section{Minimal curvature memory makes geometry a state variable}
\label{sec:memory}

A minimization algorithm that computes every new search direction from the
current gradient alone has no explicit record of the path already followed.
In a narrow valley tilted relative to the coordinate axes, successive
gradient steps can therefore alternate across the valley, repeatedly
correcting motions made in earlier iterations. One way to reduce this
zigzagging is to retain information about previous search directions. We
consider the simplest form of such memory: keeping only the direction used
in the preceding accepted step.

We label successive steps of the minimization by an integer \(t\). At step \(t\), the current position is \(x_t\), the gradient is \(g_t=\nabla f(x_t)\), and the Hessian is \(H_t=H(x_t)\). The subscript simply indicates that these quantities are evaluated at the current point.

We now choose a value \(h_t\) of the exponent introduced in Section II. The corresponding search direction, obtained directly from the current gradient, has components

\[
q_{t,i}=-\operatorname{sign}(g_{t,i})\,|g_{t,i}|^{h_t}
\]

We call this the fresh direction, meaning that it is calculated using the current gradient alone, before taking account of any information retained from previous steps. The direction actually used for the next
step, denoted by $p_t$, is obtained by correcting $q_t$ with the previously
accepted search direction $p_{t-1}$. In particular,
$x_t=x_{t-1}+\alpha_{t-1}p_{t-1}$, where $\alpha_{t-1}$ is the preceding
line-search step length.
The reason for making this correction is easiest to understand for a
quadratic cost function. After an exact line minimization along
$p_{t-1}$, the cost is stationary in that direction: moving a little
farther forward or backward along the same line would give no first-order
improvement; we would like to avoid undoing that part of the minimization.
For a quadratic function with Hessian $H$, this happens if the new
direction is \emph{conjugate} to the previous one, meaning that
$p_t^THp_{t-1}=0$. Conjugacy is a curvature-weighted form of
orthogonality; it differs from ordinary perpendicularity because the
Hessian accounts for the different curvatures and couplings between
coordinates.

We impose the corresponding condition using the Hessian at the current
point. Starting from $q_t$, we add a multiple of the remembered direction:
\begin{equation}
 p_t=q_t+\beta_t p_{t-1},\qquad
 \beta_t=-\frac{q_t^TH_tp_{t-1}}{p_{t-1}^TH_tp_{t-1}}.
 \label{eq:memory}
\end{equation}
We apply the correction only if the denominator exceeds $10^{-20}$ in the numerical experiments. When the correction is applied, we have
\begin{equation}
 p_t^TH_tp_{t-1}=0.
 \label{eq:conj}
\end{equation}
For a fixed quadratic Hessian, this condition ensures that a subsequent
step along $p_t$ does not reintroduce a gradient component along
$p_{t-1}$. Indeed, $g_{t+1}=g_t+\alpha_t Hp_t$ implies
$g_{t+1}^Tp_{t-1}=0$ whenever the preceding line search was exact.
For nonlinear functions the Hessian changes along the step, so conjugacy at the starting point does not in general preserve stationarity along the previous direction.

We can also check that the corrected direction is a descent direction. An exact line search
along the preceding direction gives $g_t^Tp_{t-1}=0$. Consequently, for
the power-gradient choice of Eq.~\eqref{eq:power},
\begin{equation}
 \begin{aligned}
 g_t^Tp_t&=g_t^Tq_t\\
          &=-\sum_i |g_{t,i}|^{1+h_t}<0,
 \end{aligned}
 \label{eq:memorydescent}
\end{equation}
provided $g_t\ne0$. 
Thus, when the preceding line minimization is exact, the memory correction preserves descent in exact arithmetic. In practical implementations, however, the line search may stop before reaching the precise minimum along the chosen direction. The gradient can then retain a component along the preceding direction, and the memory correction may no longer preserve descent. Numerical roundoff and noisy gradient estimates can have similar effects.
We therefore use a safeguard: if the corrected direction is
not sufficiently downhill (see Sec.~\ref{sec:binaryoracle}), we discard the correction and use $q_t$.
The first step likewise uses $p_0=q_0$, since no preceding direction
exists.

Computing the memory correction does not require
the full Hessian matrix. We only need its product
with the remembered direction,
$u_t=H_tp_{t-1}$, from which $\beta_t$ follows
through the two scalar products $q_t^Tu_t$ and
$p_{t-1}^Tu_t$.
The same vector $u_t$ can also be used to correct
several candidate fresh directions at the same
point, since only $q_t$ changes.

The one-vector construction described above is closely related to curvature-conjugate
search directions and to methods with variable preconditioning
\citep{HestenesStiefel1952,notay2000flexible,knyazev2007variable}.
It enforces conjugacy explicitly only with the immediately
preceding direction. 
For a positive-definite quadratic cost function,
our memory correction reduces to linear conjugate
gradients if we use $h_t=1$ at every step, start
without a remembered direction, and perform exact
line minimizations. Conjugacy with all earlier
directions then follows automatically, and in exact arithmetic
the minimum is reached in at most $d$
steps~\citep{HestenesStiefel1952}.
When other exponents are introduced, this global conjugacy
property is no longer guaranteed; the memory correction
enforces only conjugacy with the immediately preceding
direction.
Here the fresh gradient
transformation controls the relative importance of gradient components,
while the memory correction responds to the coupling of directions by
local curvature. 
In particular, the next direction now depends on the remembered
$p_{t-1}$ as well as the present point. By the algorithm's \emph{state} we mean the information required to determine its next update. Once $h_t$ is given, the following state suffices for the direction update:
\begin{equation}
 S_t=(x_t,p_{t-1}).
 \label{eq:state}
\end{equation}
The choice of $q_t$ thus affects the direction remembered at the next iteration.

\section{When the best immediate step is not the best first step}
\label{sec:greed}

\subsection{Choosing a direction using only the next step}

We have two ways of modifying a gradient step: we can change the relative sizes of the gradient components or correct the direction using the preceding step. A simple strategy to combine them is to try several candidate directions and select the one expected to give the largest immediate reduction in cost. This is often called a \emph{greedy} choice.

The local quadratic approximation introduced in
Sec.~\ref{sec:local} provides a convenient way to make
this comparison. For each candidate value of the exponent
$h$, we can construct the fresh direction $q_t(h)$ and
evaluate its predicted decrease using Eq.~\eqref{eq:D}.
Alternatively, we can first apply the curvature-memory
correction of Sec.~\ref{sec:memory}, obtaining the
direction $p_t(h)$ that the algorithm would actually
follow, and compare the predicted decreases of these
corrected directions. The two criteria need not select
the same exponent, because the memory correction affects
different candidate directions by different amounts.

At first sight the second criterion seems preferable, since the corrected direction is the one we actually use. We could go further still, replacing
the quadratic approximation with an actual line
minimization along every candidate direction. This
would identify the best immediate step within the
available choices. For a quadratic cost function with
positive-definite Hessian, the quadratic prediction
and the exact line-minimization result coincide.
None of the above criteria considers what happens next:
since the direction chosen for the current step becomes
the memory used to correct the following one, a choice that looks best now may turn
out to be worse after two steps.
We give an example below in which this happens for a strictly convex quadratic.

\subsection{An exact three-dimensional counterexample}

For a positive-definite quadratic cost function the line minimization is described exactly by the quadratic model, so we can compare the choices analytically.

Consider the strictly convex function
\begin{equation}
 f(x)=\frac12 x^T Hx,\qquad H=\operatorname{diag}(1,4,16).
 \label{eq:counterH}
\end{equation}
Its minimum is at the origin. We examine the point
\begin{equation}
 x_0=x(a)=\left(a,\frac14,1\right),\qquad
 g_0=(a,1,16),
\end{equation}
with $a>0$. We imagine that the algorithm has just completed
a line minimization along the direction
\begin{equation}
 r=(0,16,-1),
\end{equation}
which it now retains in memory. This is consistent with the
preceding step: $g_0^Tr=16-16=0$, precisely the stationarity
condition following an exact line search. For example, $x_0$
can be reached from $x_0-r$ by minimizing along $r$. We label
the present point by $x_0$ only to count the two steps that
follow; the remembered direction comes from an earlier step.

At $x_0$ we compare two choices for the fresh direction:
GD and SIGN. They give
\begin{equation}
 q_G=-g_0,\qquad q_S=-\sign(g_0)=(-1,-1,-1).
\end{equation}
For either choice, the direction actually followed is obtained
by applying the same memory correction,
Eq.~\eqref{eq:memory}, with $p_{t-1}=r$.
Thus, a ``GD-first'' or ``SIGN-first'' trajectory refers
to the chosen fresh direction; the first step is taken
along its memory-corrected counterpart.

At $a=3/2$, the exact quadratic decreases are
\begin{equation}
\begin{array}{lcc}
 & \GD & \SGN\\
\text{Fresh }D(q) & 8.1919 & 8.1488\\
\text{Corrected }D(p) & 8.2633 & 9.1098
\end{array}
\end{equation}
The fresh directions slightly favor GD, whereas the
directions actually used after the memory correction favor
SIGN. The latter comparison is also an exact comparison
of the improvements obtainable by line minimization,
because the cost function is quadratic. Starting from
$f_0=9.25$, the first steps give
\begin{equation}
 f_1^{G}=0.986685,\qquad f_1^{S}=0.140152.
\end{equation}
SIGN has clearly done better on the first step.

We now follow both branches for one further step. In each
branch we recompute the gradient and compare the GD and
SIGN candidates, each corrected against the direction taken
in the first step. We then choose whichever gives the lower
cost after an exact line search. There are thus four possible
two-step sequences, $GG$, $GS$, $SG$, and $SS$, where the
letters refer to the successive fresh geometries. At
$a=3/2$, the costs after the second line search are
\begin{equation}
\begin{array}{c|cc}
 & G & S\\ \hline
 G & 0.0011701 & 0.63859\\
 S & 0.0049086 & 0.05173
\end{array}
\end{equation}
where the first action labels the rows and the second labels
the columns. Thus GD is the better second choice in either
branch, and the ratio of costs is
\begin{equation}
 \frac{f_{SG}}{f_{GG}}\simeq4.20.
 \label{eq:42}
\end{equation}
Although SIGN produced the better first step, the lowest
cost after two steps is obtained by starting with GD.
The second action has been optimized independently for
each branch over the same two available geometries.

The reversal is due to the different states left by the first step. Different first steps produce
both different positions and different remembered
directions. The second step must then be corrected to be
conjugate to the particular first direction used in that
branch. A large immediate improvement can therefore leave
a less favorable opportunity for the next step, even with
fixed curvature and exact line searches. The counterexample
establishes this effect for the specified memory rule and
the available GD/SIGN choices.

\begin{figure*}[t]
 \centering
 \includegraphics[width=0.99\textwidth]{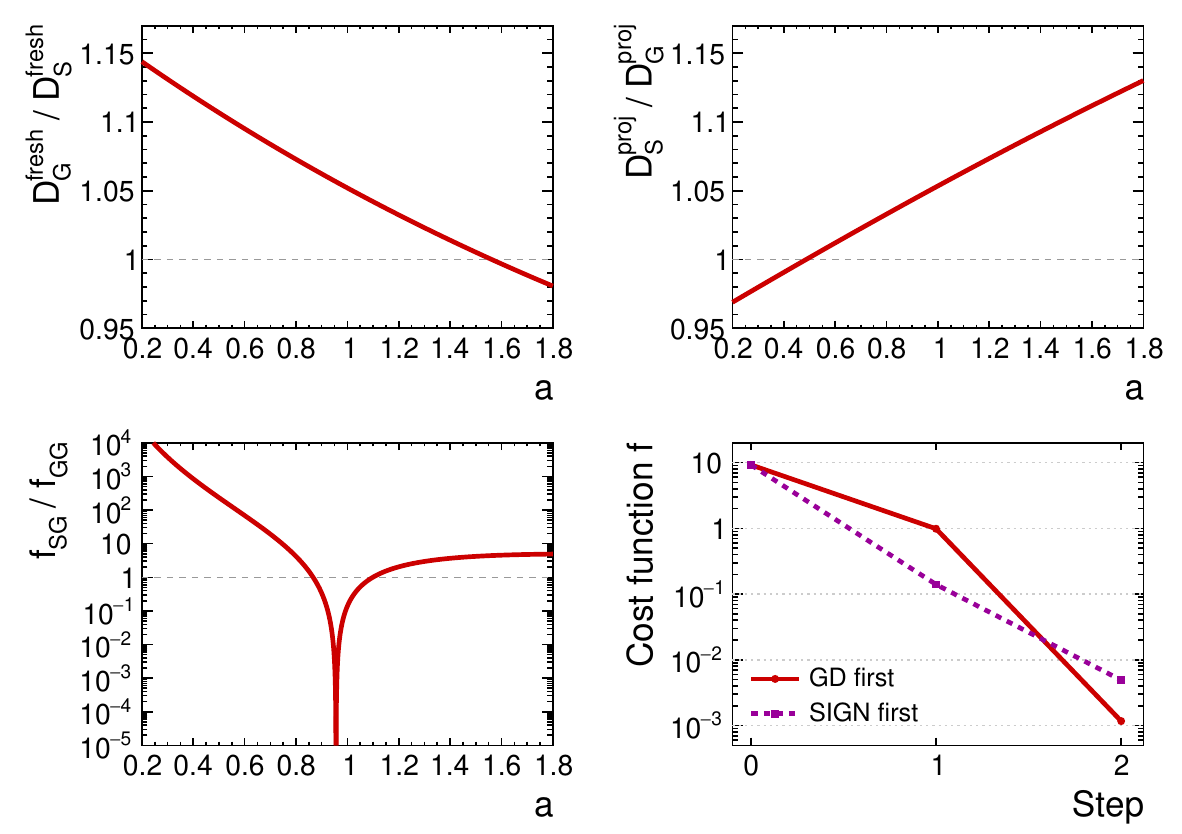}
\caption{
Comparison of immediate and two-step outcomes for
the quadratic cost function
$f(x)=\frac12 x^THx$, with
$H=\operatorname{diag}(1,4,16)$,
$x_0=(a,1/4,1)$, and remembered direction
$r=(0,16,-1)$.
The two candidate fresh directions are GD ($G$)
and SIGN ($S$), both subject to the same
curvature-memory correction.
\textit{Upper left:}
Ratio of the immediate decreases predicted for
the fresh directions, $D(q_G)/D(q_S)$;
values above unity favor GD.
\textit{Upper right:}
Ratio after memory correction,
$D(p_S)/D(p_G)$;
values above unity favor SIGN.
\textit{Lower left:}
Ratio $f_{SG}/f_{GG}$ of the costs remaining
after two steps, with the better second direction
chosen independently for each branch.
Values above unity favor starting with GD.
\textit{Lower right:}
Costs after successive steps for $a=1.5$,
comparing GD-first and SIGN-first trajectories,
each with its better second direction.
}
\label{fig:counterexample}
\end{figure*}

Figure~\ref{fig:counterexample} shows the comparison as a function of $a$ and the two trajectories for $a=3/2$. The reversal occurs over a range of values of $a$. 
The inequalities are strict at $a=3/2$, and the
reversal is not confined to this value.
Over the range $0.1\leq a\leq3$, the memory-corrected
one-step criterion favors GD below $a\simeq0.487$
and SIGN above it. The best two-step sequence,
however, starts with GD for $a<0.867$ and $a>1.099$,
and with SIGN between these values.
Thus the corrected one-step criterion makes the
wrong first choice for $0.487<a<0.867$ and
$1.099<a\leq3$.

The fresh one-step criterion favors GD up to
$a\simeq1.562$ and SIGN above this value.
It therefore disagrees with the two-step result
for $0.867<a<1.099$ and $1.562<a\leq3$.
Neither immediate criterion identifies the better
first action throughout the range examined.
The results are illustrated in Fig.~\ref{fig:counterexample}
and the corresponding expressions are given in
Appendix~A.



\section{Testing adaptive descent on nonlinear landscapes}
\label{sec:shallow}

On a quadratic cost function, the fixed Hessian allows a sequence of conjugate directions to be especially effective. In many physical minimization problems the cost is approximately quadratic only near its minimum. Farther away, the curvature may be much weaker and can vary considerably during minimization. We use three simple nonlinear examples to see whether adapting the gradient direction is useful in these circumstances.

\subsection{Three controlled families of cost functions}

Many fitting problems are handled by assigning a penalty to each measurement residual and adding these contributions to construct the total cost. Also, a residual may depend on several fit parameters; this couples their motion during minimization. Our model has both features. Denoting the benchmark cost by $F$, we write
\begin{equation}
 F(x)=\sum_{k=1}^{K}w_k\rho_k^2\,
 \ell\!\left(\frac{a_k^Tx}{\rho_k}\right).
 \label{eq:ridge}
\end{equation}
Here $x$ contains the parameters to be minimized, $a_k$ specifies a
linear combination of those parameters, $w_k>0$ gives its relative
weight, and $\rho_k>0$ sets the scale at which its contribution changes
character. The dimensionless quantity $u_k=a_k^Tx/\rho_k$ plays the
role of a simplified residual, and the function $\ell(u_k)$ determines
its penalty. Such a sum is often called a \emph{ridge-function model},
because each term depends on one projection of the parameter vector.
We use six parameters ($d=6$) and twelve randomly
oriented projections ($K=12$); their dense combinations
produce coupling among the coordinates.
The directions $a_k$ are sampled uniformly on the unit
sphere, while the weights $w_k$ and scales $\rho_k$
are sampled independently from log-uniform distributions
over $[0.3,3]$ and $[0.15,2.5]$, respectively.

We compare three choices of $\ell$, that share the same behavior for small $u$, namely $\ell(u)=u^2/2+\mathcal O(u^4)$. Consequently, all three
costs have a quadratic region around their shared minimum at $x=0$,
\begin{equation}
 F(x)=\frac12x^TH(0)x+\mathcal O(\|x\|^4),\qquad
 H(0)=\sum_k w_k a_k a_k^T,
 \label{eq:ridgecore}
\end{equation}
where $H(x)=\nabla^2F(x)$ is the Hessian introduced in
Sec.~\ref{sec:local}. For the same directions, weights, and scales,
the Hessian at the minimum is identical across the three families.
Their differences appear as we move away from that minimum.

\paragraph{Pseudo-Huber.}
The first penalty we study is
\begin{equation}
 \ell_{\rm PH}(u)=\sqrt{1+u^2}-1.
 \label{eq:ph}
\end{equation}
For small residuals this function behaves like the usual squared penalty, while for
large $|u|$ it grows approximately linearly, so its slope approaches a
constant magnitude and its curvature decreases as $|u|^{-3}$.
Large residuals therefore exert less influence than they would in a
quadratic least-squares penalty. A penalty of this form appears, up
to scaling and additive constants, in the  track-reconstruction
chain of the ANTARES experiment~\cite{ANTARES2012PointSources}, a neutrino telescope deployed in the Mediterranean sea.

\paragraph{Subquadratic, $p=1.5$.}
Our second choice is
\begin{equation}
 \ell_p(u)=\frac{(1+u^2)^{p/2}-1}{p},\qquad p=1.5.
 \label{eq:p15}
\end{equation}
For large residuals the penalty grows as $|u|^{3/2}$. Its slope in absolute value continues to increase as $|u|^{1/2}$, unlike that of pseudo-Huber, while the curvature falls as $|u|^{-1/2}$.  Related subquadratic penalties have applications in geophysical inverse problems. For example, an $L_{1.5}$ regularization term has been used to reconstruct the electrical conductivity of the Earth's mantle from geomagnetic observatory measurements, balancing the preservation of sharp conductivity contrasts against excessive smoothing \cite{Munch2018Geomagnetic}. Our smooth penalty retains the same large-amplitude power-law behavior while remaining quadratic near the origin.
Including this penalty lets us check whether the usefulness of adaptation depends on saturation of the gradient contributions.

\paragraph{Log-cosh.}
Finally, we use
\begin{equation}
 \ell_{\log\cosh}(u)=\log\cosh u.
 \label{eq:logcosh}
\end{equation}
This penalty term also becomes approximately linear for large residuals,
but its curvature decreases exponentially, rather than as a power law. A log-cosh penalty has been used for energy reconstruction in
neutrino-telescope studies~\cite{Yu2023TriggerReconstruction}. 

The three families above let us change the behavior of the cost far from its minimum while retaining the same quadratic core and the same
underlying parameter couplings.
Our benchmarks are intentionally simpler than actual fitting problems. Linearizing a measured residual near a reference parameter set gives
$r_k(\theta_0+\Delta\theta)\simeq b_k+a_k^T\Delta\theta$.
Our examples take the offsets $b_k$ to be zero, so that every term has
its minimum at the same known point. Actual residuals may have nonzero offsets and depend nonlinearly on the parameters; here we concentrate on the effects of changing curvature.

We also need a way to compare starting points in different random models. The characteristic scales of the models vary, so the same distance from
the origin need not imply a comparable degree of nonlinearity. We
instead monitor the total local curvature through the trace of the
Hessian, which is the sum of its diagonal entries. We define
\begin{equation}
 \Phi(x)=\frac{\tr H(x)}{\tr H(0)}.
 \label{eq:phi}
\end{equation}
In the three families, $\Phi$ approaches one near the quadratic
minimum and becomes smaller as the average curvature weakens.
Here \emph{shallow} refers to this reduction in curvature; it does
not imply that the gradient itself must be small.

For each randomly generated cost function, we sample
a unit direction $u$ uniformly on the sphere and
determine the starting point $x_0=ru$ by adjusting
the radius $r$ through bisection until the prescribed
value of $\Phi_0$ is attained.
Repeating this procedure for different values of \(\Phi_0\) allows us to compare the behavior of the minimization algorithms across curvature regimes.
The quantity $\Phi$ is available here because the minimum is known
by construction. We use it to prepare and classify the numerical
tests; a practical optimizer would generally have no access to the
denominator $\tr H(0)$.


\subsection{An idealized two-step lookahead}
\label{sec:binaryoracle}

Testing a particular rule for switching between GD and SIGN tells us little about how much could have been gained with better choices. To get an idea of that gain, we try several short sequences of steps from the same state and compare the final costs before choosing the next step. We call this procedure an \emph{oracle}, since it uses information that would normally be expensive to obtain during minimization.

We first restrict the available fresh directions to GD and SIGN,
corresponding to the exponents $h=1$ and $h=0$ of Eq.~\eqref{eq:power}. The curvature-memory correction of
Sec.~\ref{sec:memory} and an exact line minimization are applied
at every step, irrespective of the chosen fresh direction. We
compare three methods: GD-memory, which always proposes GD;
SIGN-memory, which always proposes SIGN; and an oracle that is
allowed to alternate between the two. Thus ``fixed GD'' and
``fixed SIGN'' describe the proposed directions: both methods still
use curvature memory.

To decide the next move at iteration $t$,
the oracle examines the consequences of two successive steps. With GD and SIGN as the available choices, there are four possible sequences of fresh directions:

\[
GG,\qquad GS,\qquad SG,\qquad SS.
\]

The first letter specifies the direction proposed for the next step; the second specifies the direction proposed for the step after that.
Each trial starts from the same current point and the
same remembered direction. Within each trial, both steps receive
their own memory corrections and exact line minimizations. The
oracle compares the true cost reached after the two trial steps,
and chooses the sequence giving the lowest value.

The algorithm then executes only the first step of the preferred sequence. At the new point, it evaluates the four
possibilities again, using the updated gradient, curvature, and
remembered direction. The second action originally planned may
therefore differ from the one actually taken. This is a \emph{receding-horizon} procedure. It can finish worse than a fixed-geometry method because it does not optimize the entire six-step trajectory at once.

For each randomly generated landscape and prescribed starting value
$\Phi_0$, we run the oracle and the two fixed-geometry methods from
exactly the same initial point. All runs begin without a remembered
direction; the first step therefore uses the fresh direction
unchanged. We compare the costs after six accepted steps, denoting
them by $F_6^{\GD}$, $F_6^{\SGN}$, and $F_6^{\mathrm{or}}$.
We compare the methods at a fixed number of accepted steps. The oracle is clearly more expensive, as it explores four speculative two-step trajectories at each decision. We use it to assess the potential gains from adaptation.

Each line minimization brackets the minimum by
successively doubling the step length, followed by
72 bisection iterations based on the sign of the
directional derivative. The memory correction is
discarded, reverting to $p_t=q_t$, whenever
\[
g_t^Tp_t\geq-10^{-12}\max(1,\|g_t\|\,\|p_t\|).
\]
When candidate sequences yield numerically identical
final costs, the oracle selects the first among
them in the ordering $GG$, $GS$, $SG$, $SS$.
 
We use a \emph{paired} comparison: each method is tested on the same
random problem and starting point. For each case, we compare the
oracle with whichever of the two fixed-geometry trajectories finishes
with the lower cost. The reported advantage at a given $\Phi_0$ is
\begin{equation}
 A(\Phi_0)=\operatorname{median}_{j}\!\left[
 \frac{\min\!\left(F_{6,j}^{\GD},F_{6,j}^{\SGN}\right)}
      {F_{6,j}^{\mathrm{or}}}\right],
 \label{eq:advantage}
\end{equation}
where $j$ labels independently generated random problems. The
median is taken \emph{after} forming the ratio for each matched
case. Values above one indicate that the oracle gives a smaller
final cost on a typical case, even when compared with the better
fixed strategy for that case. Because that reference is chosen
after observing both fixed-method outcomes, it is a deliberately
demanding benchmark. We also report the fraction of cases
for which the oracle finishes strictly below \emph{both} fixed
strategies. This second measure indicates how frequently an
improvement occurs, as distinct from its typical relative size.

We study 23 prescribed values of $\Phi_0$ in each landscape family.
The pseudo-Huber scan contains 1000 independent random problems at
each value, while the subquadratic and log-cosh scans contain 3000
each. The random constructions are matched across the families,
with the first 1000 problems common to all three, and the same
random model and starting ray are reused across the prescribed
values of $\Phi_0$. Uncertainties on the reported medians are
estimated with 68\% bootstrap confidence intervals (obtained by
resampling the random problems); uncertainties on the paired win
fractions use 95\% Wilson binomial confidence intervals.


\subsection{From shallow landscapes to quadratic minima}
\label{sec:phase}

\begin{figure*}[t]
 \centering
 \includegraphics[width=0.99\textwidth]{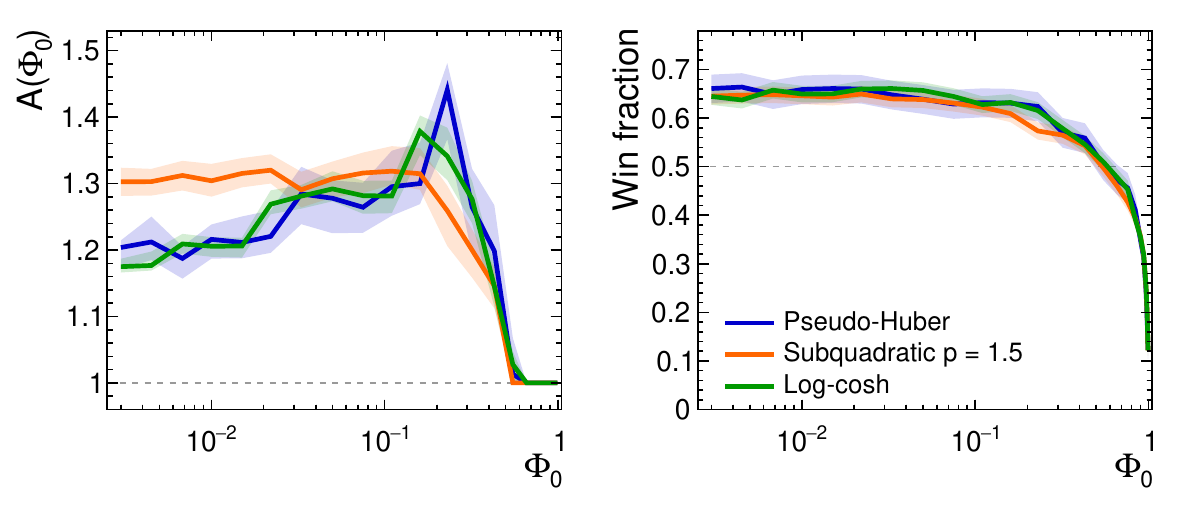}
\caption{
Performance of the binary two-step oracle after six
accepted steps, as a function of the initial curvature ratio $\Phi_0=\operatorname{tr}H(x_0)/\operatorname{tr}H(0)$.
Smaller values of $\Phi_0$ correspond to shallower regions.
\textit{Left:}
Median paired advantage $A(\Phi_0)$, defined in
Eq.~\eqref{eq:advantage}, relative to the better
of GD-memory and SIGN-memory on each problem.
\textit{Right:}
Fraction of problems on which the oracle finishes
strictly below both fixed methods.
The horizontal dashed line marks $50\%$.
Shaded bands show $68\%$ bootstrap confidence
intervals for the median advantage and $95\%$
Wilson intervals for the win fraction.
The pseudo-Huber ensemble contains 1000 problems
per $\Phi_0$; the subquadratic and log-cosh
ensembles contain 3000 each.
} \label{fig:phase}
\end{figure*}

The results of the two-step lookahead experiment are summarized in
Table~\ref{tab:phase} and Fig.~\ref{fig:phase}. 
The left panel shows the median paired advantage
$A(\Phi_0)$ defined in Eq.~\eqref{eq:advantage},
while the right panel shows the fraction of
problems on which the oracle outperforms both
fixed-geometry methods. The two panels tell us
how large the typical advantage is and how
often an improvement occurs.

In the shallow part of the landscapes, where $\Phi_0$ is small, both
measurements favor adaptation. The oracle beats the better fixed
strategy on roughly $60$--$66\%$ of the random problems over much of
this region. The size of the improvement depends more strongly on the
form of the cost function. The largest median advantage is about
$1.44$ for pseudo-Huber, $1.32$ for the subquadratic family, and
$1.38$ for log-cosh. A value $A=1.44$, for example, means that the
better fixed strategy leaves a final cost about $44\%$ higher than the
oracle on the median paired comparison. The advantage has a pronounced
maximum at an intermediate $\Phi_0$ for pseudo-Huber and log-cosh,
whereas the subquadratic family exhibits a broader plateau toward the
shallow end of the scan.

\begin{table*}[t]
\centering
\caption{Binary two-step oracle results after six accepted
steps. The maximum median paired advantage $A$, its location, and the
value $\Phi_c$ at which the paired win fraction crosses $50\%$ are
obtained from the fine scans. The crossover positions are estimated
by linear interpolation between adjacent sampled values of $\Phi_0$;
they are descriptive markers for these ensembles and settings. The
last column gives the win fraction at the sampled point closest to
the quadratic minimum, $\Phi_0=0.995$.}
\label{tab:phase}
\small
\begin{tabular*}{0.95\textwidth}
  {@{\extracolsep{\fill}}lcccc@{}}
\toprule
Landscape
 & \shortstack{Max.\\$A$}
 & \shortstack{$\Phi_0$ at\\max. $A$}
 & \shortstack{$\Phi_c$\\(50\% wins)}
 & \shortstack{Win fraction\\at $\Phi_0=0.995$} \\
\midrule
Pseudo-Huber
 & 1.441 & 0.230 & 0.561 & 0.123 \\
Subquadratic $p=1.5$
 & 1.320 & 0.022 & 0.543 & 0.122 \\
Log-cosh
 & 1.379 & 0.160 & 0.584 & 0.122 \\
\bottomrule
\end{tabular*}
\end{table*}

As the starting point approaches the quadratic region, the benefit
becomes less frequent. The fraction of problems on which the oracle
beats \emph{both} fixed methods falls through $50\%$ at
$\Phi_0\simeq0.54$--$0.58$ in all three families. This close agreement
is notable: the three penalties have different large-residual
behavior, including a subquadratic model whose gradient does not
saturate. The crossover is therefore present beyond models with
asymptotically constant gradient contributions. Also, the common random geometry and memory construction
limit how broadly these numerical crossover values
should be interpreted.
The crossing of the $50\%$ win fraction should not
be interpreted as the point at which losses begin
to outnumber wins, since ties count as non-wins.
For example, at $\Phi_0=0.55$ the pseudo-Huber
ensemble exhibits win and loss fractions of $0.50$
and $0.29$, respectively.
Only much closer to the quadratic regime, between
$\Phi_0=0.90$ and $0.94$, do losses become more
frequent than wins in all three families.

Very close to the quadratic minimum, the paired comparison requires
some care. At $\Phi_0=0.995$ the median advantage equals one in all
three families, which might suggest that choosing adaptively has
little effect. Looking at individual outcomes, we find that the oracle improves upon both fixed methods in only about $12\%$ of the
trials, gives a higher final cost than the better fixed method in
roughly $46\%$, and ties it in the remaining $41$--$42\%$. The median
ratio equals one because many trajectories reproduce the fixed result exactly; it does not imply that unfavorable interventions have disappeared. 

Curvature memory may explain part of this behavior. Near a quadratic minimum, the Hessian changes little between steps, and
maintaining curvature-conjugate directions becomes particularly
valuable; repeatedly changing the fresh gradient transformation may
interfere with that structure. Our experiments cannot isolate this mechanism, since each choice changes both the position and the direction retained in memory.
Nevertheless, the observed deterioration is consistent with the
importance of preserving a favorable sequence of
directions when the cost becomes nearly quadratic.
A further consideration in interpreting these results
is that each run consists of six accepted steps in
a six-dimensional parameter space.
On a positive-definite quadratic, GD-memory with
exact line minimizations reduces to linear conjugate
gradients, which terminates in at most six steps
in exact arithmetic.
Our choice of step budget may therefore particularly
favor GD-memory in the near-quadratic regime.
The quantitative extent of the observed deterioration
should be interpreted in light of this experimental
setting.

Switching between GD and SIGN can thus be beneficial, particularly in the shallow regions of our examples. In the next section we examine whether intermediate gradient powers can improve on these two choices.


\section{Intermediate power gradients and changing curvature}
\label{sec:continuous}

\subsection{Are intermediate gradient powers useful?}
\label{sec:powergrid}

In Section~\ref{sec:phase} we showed that changing between GD and SIGN
can improve optimization in the shallow regions of our benchmark
landscapes. These are two extreme ways of using the gradient: GD
preserves the relative sizes of its components, while SIGN discards
them. There is no reason, however, to expect the most useful choice
to lie at either extreme. For example, if two gradient components have
magnitudes 100 and 1 in our fixed, scaled coordinates, GD gives their
directional components a ratio of 100 to 1, whereas SIGN gives them
equal magnitudes. A power exponent $h=1/2$ produces the intermediate
ratio 10 to 1, reducing the dominance of the larger gradient component without treating both components identically.

To test whether such intermediate choices can improve a trajectory,
we extend the candidate exponents to
\begin{equation}
 \mathcal E=\{0,\,0.25,\,0.5,\,0.75,\,1\}.
 \label{eq:hgrid}
\end{equation}
We use the same curvature-memory correction and exact line
minimizations as before. This experiment differs from the receding-horizon runs of
Sec.~\ref{sec:binaryoracle}, because here we first generate a trajectory using
GD-memory alone; then, after each of its first five accepted steps, we
\emph{freeze} the state (the current position and the direction retained from the preceding step). From that identical state we then
try every possible pair of exponents in $\mathcal E$, applying the
memory correction and an exact line minimization at both steps.
There are $5\times5=25$ two-step continuations. The trials never alter the underlying GD-memory trajectory. For each cost-function family we use
1000 random model realizations across the same 23 prescribed starting
curvature values introduced in Sec.~\ref{sec:shallow}.

Let $h_a$ denote the exponent used on the first trial step and
$h_b$ the exponent on the second. For a fixed first choice, we find
the lowest final cost obtainable by choosing the second exponent
from the same set:
\begin{equation}
 V_t(h_a)=\min_{h_b\in\mathcal E}
 F\!\left(x_{t+2}(h_a,h_b)\right).
 \label{eq:Vh}
\end{equation}
Here $x_{t+2}(h_a,h_b)$ is the point reached after performing both
trial steps, including their memory corrections. Thus $V_t(h_a)$ is the lowest cost after two steps for a given first exponent. We denote the preferred first exponent within this
restricted two-step search by
\begin{equation}
 h_*^{(2)}\in\argmin_{h_a\in\mathcal E}V_t(h_a).
 \label{eq:hstar}
\end{equation}
The superscript indicates the two-step comparison.
When two or more first exponents give final costs
within $10^{-10}F(x_t)$ of the minimum, we treat
them as tied rather than selecting a winner
according to numerical roundoff. Each tied
exponent receives equal statistical weight.
We denote the arithmetic mean of the tied exponents at a given state by $\bar h_t^{(2)}$. It coincides with $h_*^{(2)}$ when only one exponent lies within the tie tolerance. We use $\bar h_t^{(2)}$ when computing averages over states.

We also need to distinguish a genuinely useful interior first
choice from a case in which only the second step benefits from
an intermediate exponent. For this purpose, we compare the best
five-level continuation with the best sequence whose \emph{first}
exponent is restricted to GD or SIGN. Both are allowed to choose
any of the five exponents on their second step. We define the
additional benefit of allowing an interior first exponent as
\begin{equation}
 \Delta_{\mathrm{int}}(t)=
 \frac{\min\{V_t(0),V_t(1)\}
       -\min_{h\in\mathcal E}V_t(h)}{F(x_t)}.
 \label{eq:intgain}
\end{equation}
The numerator is the further reduction in the final two-step cost
made possible by extending the \emph{first} choice beyond the two
endpoints. Dividing by the current cost makes the gain dimensionless
and allows it to be compared across different random problems.
By construction, $\Delta_{\mathrm{int}}\geq0$.

\begin{table*}[t]
\centering
\caption{Preferred first exponent in the restricted five-level,
two-step search, evaluated at frozen GD-memory states. The columns
labelled ``shallow'' and ``near quadratic'' use the \emph{current}
curvature ratio, $\Phi(x_t)<0.55$ and $\Phi(x_t)\geq0.84$,
respectively; these are descriptive subdivisions of the sampled
states. ``Interior first'' is the tie-weighted fraction of states
preferring $h=0.25$, $0.5$, or $0.75$. The final column is the mean
additional benefit of allowing an interior first exponent, as
defined in Eq.~\eqref{eq:intgain}, including states where that
benefit vanishes. The study uses 1000 random models at each
prescribed starting curvature value, following five GD-memory
steps per starting trajectory.}
\label{tab:hstruct}
\small
\begin{tabular}{lccccc}
\toprule
Landscape & $\langle\bar h_t^{(2)}\rangle$ all &
 $\langle\bar h_t^{(2)}\rangle$ shallow &
 $\langle\bar h_t^{(2)}\rangle$ near quadratic &
 Interior first &
 $\langle\Delta_{\mathrm{int}}\rangle$ shallow \\
\midrule
Pseudo-Huber & 0.661 & 0.580 & 0.737 & 35.7\% & 1.495\% \\
Subquadratic $p=1.5$ & 0.676 & 0.617 & 0.748 & 37.6\% & 1.504\% \\
Log-cosh & 0.661 & 0.576 & 0.737 & 35.6\% & 1.456\% \\
\bottomrule
\end{tabular}
\end{table*}

\begin{figure*}[ht]
 \centering
 \includegraphics[width=0.99\textwidth]{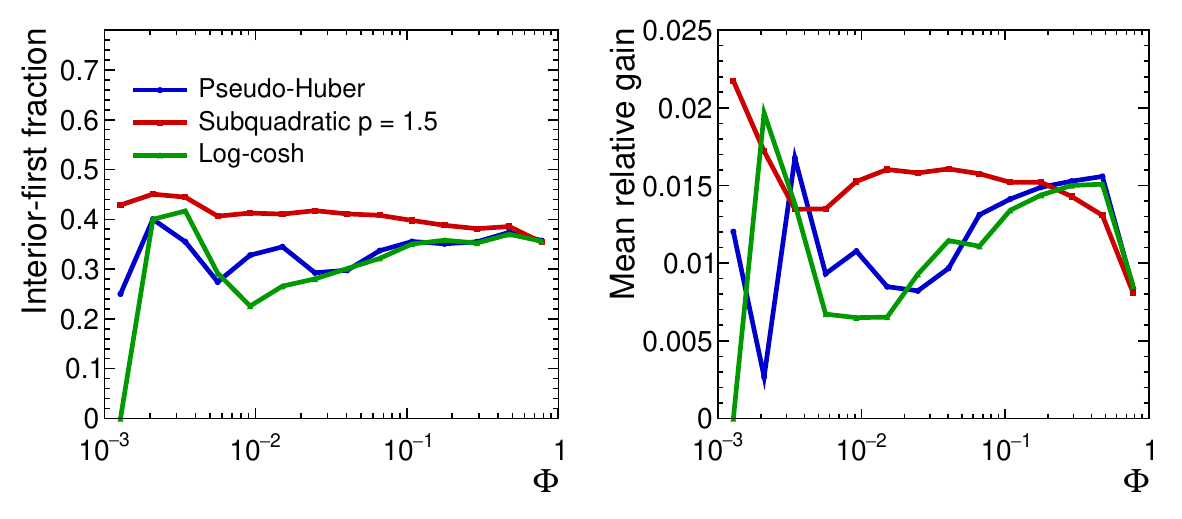}
\caption{
Value of interior gradient exponents in the
five-level two-step oracle, evaluated at frozen
GD-memory states across the three nonlinear
cost-function families.
\textit{Left:}
Tie-weighted fraction of optimal first actions
with an interior exponent $0<h<1$.
\textit{Right:}
Mean additional gain $\Delta_{\mathrm{int}}$ from allowing interior first exponents, as defined in Eq.~\eqref{eq:intgain}.
Both alternatives use the same five-level
optimization of the second action.
The horizontal axis is the curvature ratio
$\Phi$ at the frozen decision state.
}
\label{fig:interiorvalue}
\end{figure*}

The results are summarized in Table~\ref{tab:hstruct} and
Fig.~\ref{fig:interiorvalue}. 
The tie-weighted fraction favoring an interior first exponent is about $36$--$38\%$ in all three families. In the shallow region, allowing intermediate first exponents produces an additional two-step decrease averaging about 1.5$\%$ of the current cost. This average includes all sampled states in that region, including those where one of the two endpoint exponents already gives the best result. The improvement varies across individual states, but demonstrates that intermediate powers can offer a measurable advantage over restricting the first choice to GD or SIGN.

The preferred exponent also changes with the local curvature regime.
On average, smaller exponents are favored in shallower regions,
where compressing the relative gradient magnitudes can be useful.
Closer to the quadratic core, the preferred first exponent shifts
toward GD. At a given $\Phi$, though, the distribution of preferred exponents is still broad.

Intermediate powers can therefore give an advantage
over both GD and SIGN, although identifying the
best choice through a two-step search requires
25 trial continuations at each state.


\subsection{Measuring curvature changes along the preceding step}
\label{sec:nsec}

The five-level search reveals a connection between the preferred
exponent and the curvature regime, but the quantity $\Phi$ used to
classify our test problems is not generally available during an
optimization. Its definition requires the Hessian at the minimum,
which is known by construction in our benchmarks. We would like to recognize changes in curvature using quantities available during minimization.

Consider what happens during one accepted step, from $x_{t-1}$ to
$x_t$. We denote the displacement and the corresponding change in
gradient by
\begin{equation}
 s_t=x_t-x_{t-1},\qquad y_t=g_t-g_{t-1}.
 \label{eq:secsy}
\end{equation}
The gradient tells us how the cost changes with the parameters. Its
change over a step therefore contains information about the curvature
encountered along that step. If the cost were exactly quadratic, its
Hessian would be constant, and we would have $y_t=H_ts_t$ exactly.
For a general cost function, $H_ts_t$ instead predicts the gradient
change using only the curvature measured at the end of the step. Comparing these two vectors is a way of testing how well a
constant-curvature description represents the motion just completed.

We quantify the disagreement between the two vectors by the dimensionless number
\begin{equation}
 N_{\mathrm{sec}}(t)=
 \frac{\|H_ts_t-y_t\|}{\|H_ts_t\|+\|y_t\|}.
 \label{eq:nsec}
\end{equation}
The denominator supplies a scale for the comparison, so the result
does not depend on multiplying the entire cost function by a
constant. Whenever the denominator is nonzero, the triangle
inequality gives $0\leq N_{\mathrm{sec}}\leq1$; we set it to zero
if both vectors vanish. A value near zero means that the observed
gradient change agrees closely with the endpoint-Hessian prediction.
A larger value indicates a stronger disagreement, relative to the
sizes of the two vectors.
We can express the gradient change as an integral of the Hessian along the step:
\begin{equation}
 y_t=\int_0^1
 H(x_{t-1}+\tau s_t)s_t\,d\tau.
 \label{eq:secint}
\end{equation}
Subtracting this expression from $H_ts_t$ gives
\begin{equation}
 H_ts_t-y_t=\int_0^1
 \left[H_t-H(x_{t-1}+\tau s_t)\right]s_t\,d\tau.
 \label{eq:secdiff}
\end{equation}
Thus $N_{\mathrm{sec}}=0$ for any quadratic cost function, regardless
of how strongly the parameters are coupled or how different its
curvatures are along different directions. For a nonlinear cost,
the statistic measures the accumulated discrepancy between endpoint
curvature and the curvature sampled along the preceding displacement.
The diagnostic can miss variations in directions not probed by the step, or variations that cancel in the integral.

In our algorithm $N_{\mathrm{sec}}$ is inexpensive to calculate. Since
the preceding displacement is parallel to the previous search
direction, $s_t=\alpha_{t-1}p_{t-1}$, the vector $H_ts_t$ is just
$\alpha_{t-1}H_tp_{t-1}$. The Hessian--vector product
$H_tp_{t-1}$ is already required for the curvature-memory
correction in Sec.~\ref{sec:memory}. 
The gradient difference $y_t$
requires retaining the preceding gradient, but no additional
Hessian--vector product is needed. Unlike $\Phi$, this
statistic also requires no knowledge of the location or
curvature of the minimum.

The state $S_t$ defined in
Eq.~\eqref{eq:state} is sufficient for constructing
the corrected direction once the exponent is known.
A policy using $N_{\mathrm{sec}}$ must additionally
retain the preceding gradient and step length.

\begin{figure*}[t]
 \centering
 \includegraphics[width=0.99\textwidth]{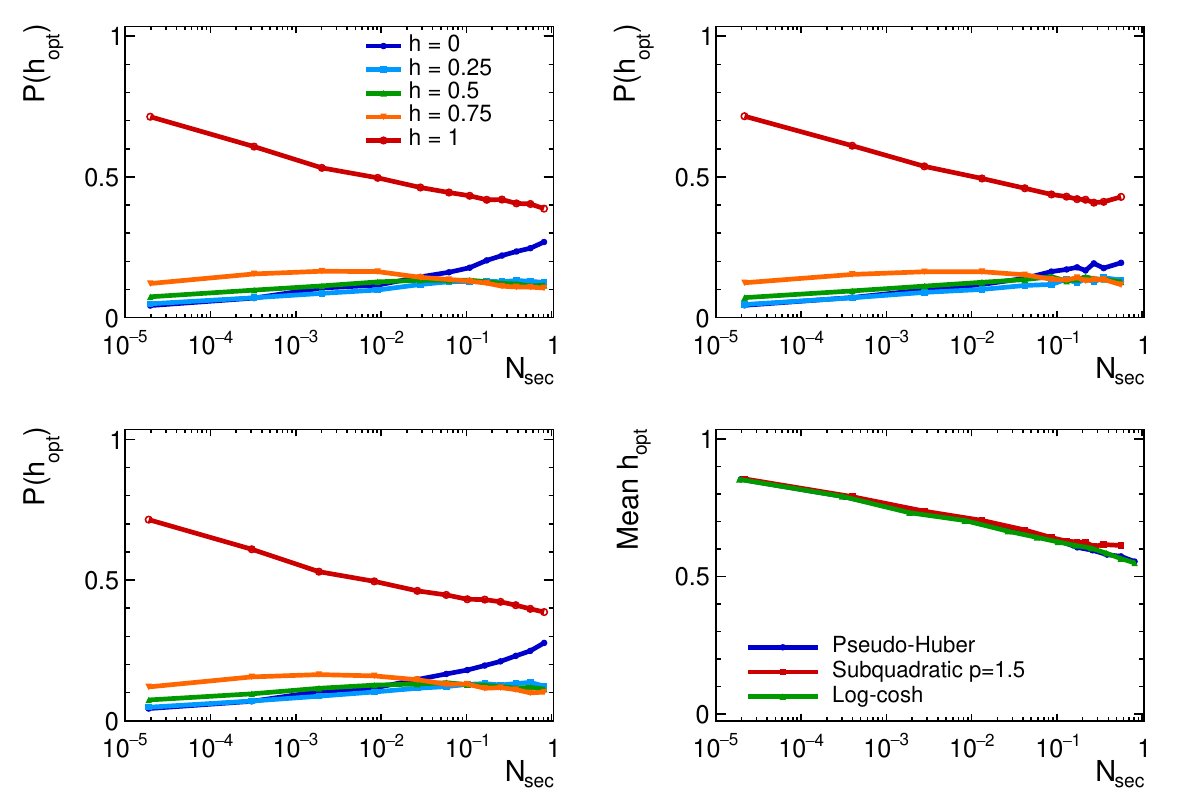}
\caption{
Dependence of the preferred first gradient exponent
on the secant nonquadraticity diagnostic
$N_{\mathrm{sec}}$, evaluated at frozen GD-memory
states using the five-level two-step oracle.
\textit{Upper left:}
Tie-weighted distribution of the optimal first
exponent for the pseudo-Huber family.
\textit{Upper right:}
The corresponding distribution for the
subquadratic family ($p=1.5$).
\textit{Lower left:}
The corresponding distribution for the
log-cosh family.
\textit{Lower right:}
Mean tie-aware optimal first exponent for
all three families.
The horizontal coordinate is the median
$N_{\mathrm{sec}}$ within each family-specific
quantile bin.
}
 \label{fig:hvsnsec}
\end{figure*}

We evaluate $N_{\mathrm{sec}}$ at the same frozen GD-memory states
used in Sec.~\ref{sec:powergrid}, and compare it with the preferred
first exponent from the restricted two-step search, using $\bar h_t^{(2)}$ for the means. The results are shown in Fig.~\ref{fig:hvsnsec}. In all
three families, states with very small mismatch favor exponents
closer to GD on average. As the mismatch increases, the mean
preferred exponent decreases, so stronger compression of gradient magnitudes is preferred on average. The mean curves are notably
similar over their overlapping ranges of $N_{\mathrm{sec}}$,
despite the different large-residual behavior of the three
penalties. In the lowest mismatch bins the mean exponent is about
$0.85$ in each family; toward the high-mismatch end of the sampled
ranges it falls to roughly $0.55$--$0.61$.

Even within a narrow range of $N_{\mathrm{sec}}$, different states favor different exponents. The average is well correlated with the mismatch, but we still have considerable uncertainty about which exponent to choose at a particular state.


\section{Predictable geometry is not the same as a useful policy}
\label{sec:policy}

The similarity of the mean curves in Fig.~\ref{fig:hvsnsec} suggests
that $N_{\mathrm{sec}}$ could help us choose the exponent
during minimization. But a good prediction of the
preferred exponent does not necessarily mean that
using it will improve a complete run. The states
examined in Sec.~VI came from GD-memory trajectories,
whereas an adaptive method generates its own states.
We examine both aspects below.


\subsection{Predicting the preferred exponent across landscapes}
\label{sec:losopred}

To test whether the observed relation transfers to another family, we fit the relation between $N_{\mathrm{sec}}$ and
the preferred exponent using two families, and compare its predictions with the third. We repeat this leave-one-family-out test with each family held out in turn.

The quantities to be predicted come from the frozen-state
experiments of Sec.~\ref{sec:powergrid}. At each state, the
restricted two-step search identifies which of the five candidate
first exponents gives the smallest cost after its best available
second step. As described there, nearly equal outcomes share
statistical weight. We are interested in the \emph{average}
preferred exponent among states with comparable values of
$N_{\mathrm{sec}}$. In particular, we do not attempt here to
identify the best exponent for each individual state.

We first combine the frozen states from the two training families,
order them by $N_{\mathrm{sec}}$, and divide them into 24 groups
containing approximately equal numbers of states. 
Within each group we average the statewise
tie-aware preferred exponents defined in
Sec.~\ref{sec:powergrid}, obtaining the conditional
mean exponent for that range of $N_{\mathrm{sec}}$. We associate each mean with the median $N_{\mathrm{sec}}$ of its group. Equal-population bins give us comparable numbers of observations per bin, although densely sampled mismatch ranges still have more bins.

For this fit we choose a simple decreasing curve,
\begin{equation}
 h_{\mathrm{fit}}(N)=h_{\mathrm{low}}+
 \frac{h_{\mathrm{high}}-h_{\mathrm{low}}}
 {1+(N/N_0)^\gamma},
 \qquad N=N_{\mathrm{sec}}.
 \label{eq:hmap}
\end{equation}
Here $h_{\mathrm{high}}$ is the limiting exponent at very small
mismatch, $h_{\mathrm{low}}$ is the formal large-$N$ asymptote (outside the physical range $0\leq N\leq1$), $N_0$ sets the scale over which the decrease takes place,
and $\gamma>0$ controls its shape. 
The parameters are determined by minimizing the squared differences between the fitted curve and the mean preferred exponents in the 24 groups. Each squared difference is weighted by the number of states in the corresponding group. Since the groups contain approximately equal numbers of states, their contributions to the fit are nearly equal.
We choose this form empirically. The parameters are obtained by a coarse grid search followed by a finer search around its best point. Although the underlying two-step search uses five discrete exponents, the
fitted curve can take intermediate values between
zero and one. The fitted parameters are reported
in Table~\ref{tab:mapparams}, together with those
of the regret map introduced below.
For the mean map the lower asymptote lies below the observed bin means. It is a parameter of the empirical form, rather than a measured limiting exponent.

\begin{figure*}[t]
 \centering
 \includegraphics[page=1,width=0.99\textwidth]{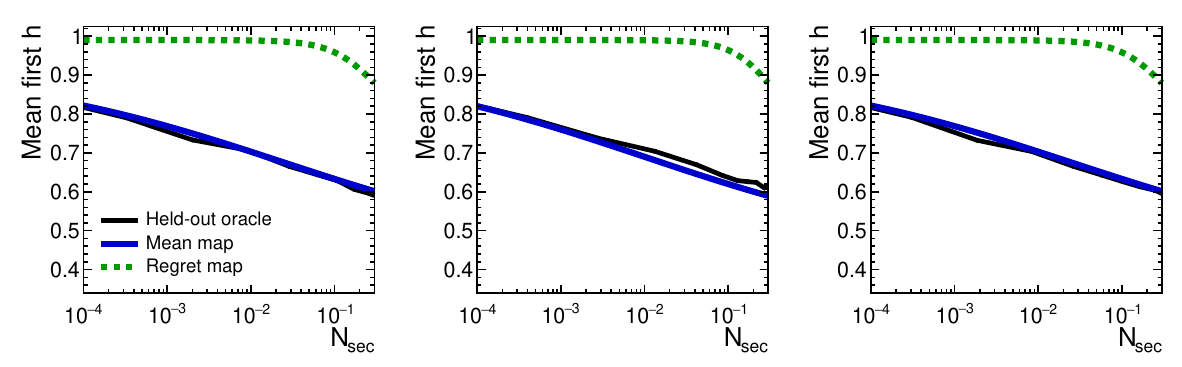}
\caption{
Leave-one-family-out transfer of the relationship
between the secant nonquadraticity diagnostic
$N_{\mathrm{sec}}$ and the preferred first
gradient exponent.
From left to right, the held-out families are
pseudo-Huber, subquadratic ($p=1.5$), and log-cosh.
The black curve shows the binned mean tie-aware preferred first exponent obtained from the two-step oracle
on the held-out family.
The blue curves are monotone sigmoid maps fitted
to the mean optimal exponent using only the
other two families.
The green dashed curves use the same functional
form but are fitted to minimize frozen-state
two-step regret.
Neither map uses data from the held-out family
during fitting.
}
 \label{fig:losostruct}
\end{figure*}

After fitting, we keep the parameters fixed and evaluate the
resulting curve on the omitted family. For this test, the states
of that family are grouped into 24 equal-population bins of
$N_{\mathrm{sec}}$, and their mean preferred exponents are
compared with the predictions at the median mismatch of each bin.
We compute the root-mean-square error (RMSE) between the binned means and predictions, weighting each bin by its population. We obtain
\begin{equation}
\begin{aligned}
 \mathrm{RMSE}_{\mathrm{PH}} &= 0.0111,\\
 \mathrm{RMSE}_{p=1.5} &= 0.0224,\\
 \mathrm{RMSE}_{\log\cosh} &= 0.0124.
\end{aligned}
\label{eq:rmse}
\end{equation}
For comparison, a constant predictor given by the
mean preferred exponent over the two training
families yields held-out RMSE values of $0.090$,
$0.079$, and $0.090$, respectively.
The fitted relation therefore reduces the prediction
error by factors ranging from approximately $3.5$
to $8$ relative to this baseline.
The three panels of
Fig.~\ref{fig:losostruct} show that the fitted curves follow the
observed mean trends even though the corresponding family was
excluded from training.

The agreement is encouraging, particularly given the different penalties in the three families. They do share the same construction of random models, though, and these small errors refer to averages over many states. As Fig.~\ref{fig:hvsnsec} shows, we cannot expect such accuracy for individual decisions.

\subsection{When a good prediction becomes a poor optimization rule}
\label{sec:closedloop}

The agreement found in Sec.~\ref{sec:losopred} suggests a direct
way to choose the exponent during minimization. At each iteration,
we measure $N_{\mathrm{sec}}$, insert it into the fitted relation
$h_{\mathrm{fit}}(N_{\mathrm{sec}})$, and use the resulting exponent
to construct the fresh direction. The curvature-memory correction and line minimization are then applied as usual. We call this the
\emph{mean-map policy}. A policy here is simply a rule for choosing an action from information available at the current step. Unlike
the oracle searches, the mean map does not test alternative future trajectories while the optimizer is running.

We test this rule on new random realizations of the held-out family, distinct from the frozen states used to determine the map.
Each method starts from the same point with a budget of six accepted steps. The first step is GD for every method, since $N_{\mathrm{sec}}$ cannot be evaluated until a step has been taken.
After the first step we evaluate the map at the current mismatch. The subsequent states are generated by this policy on new problems, rather than sampled from a GD-memory trajectory.

To measure the performance of the above rule, we compare the final cost of a method $M$
with that of GD-memory on the \emph{same} random problem. A simple
ratio of these two final costs becomes difficult to interpret when both are very close to zero; we therefore measure the signed
improvement relative to the common initial cost,
\begin{equation}
 I_{M,j}=\frac{F_{6,j}^{\GD}-F_{6,j}^{M}}{F_{0,j}},
 \label{eq:safegain}
\end{equation}
where $j$ labels a random problem and starting point. Positive values indicate an improvement over GD-memory; negative values indicate that the proposed adaptation has left a higher cost.
This normalization remains well behaved even when the final cost of GD-memory is extremely small.

For a compact comparison across all tested starting-curvature
values, we add these signed improvements and express their sum
relative to the total \emph{positive} gain available from the
five-level two-step oracle on the same test cases:
\begin{equation}
 R_M=\frac{\displaystyle\sum_j I_{M,j}}
 {\displaystyle\sum_j\max\{0,I_{\mathrm{or},j}\}},
 \qquad
 I_{\mathrm{or},j}=
 \frac{F_{6,j}^{\GD}-F_{6,j}^{\mathrm{or}}}{F_{0,j}}.
 \label{eq:saferecovery}
\end{equation}
The denominator gives a common scale for all methods on the same test sample. Both gains and losses enter the numerator, so a few large losses can outweigh many small gains. A value of zero means no net improvement over GD-memory. Even the oracle can have $R_{\mathrm{or}}<1$ because it sometimes loses to GD-memory.

\begin{figure*}[t]
 \centering
 \includegraphics[width=0.99\textwidth]{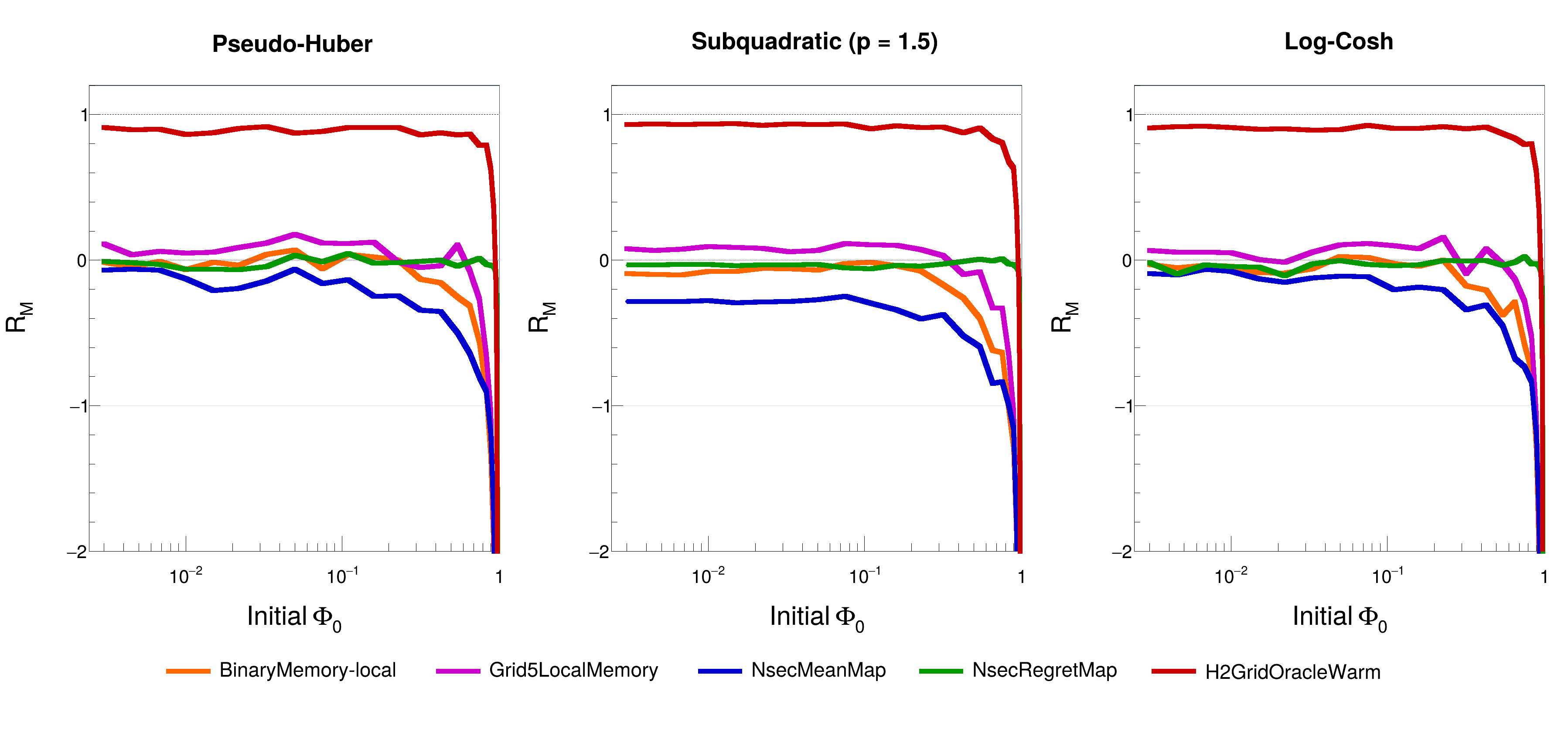}
 
\caption{Closed-loop comparison of adaptive direction
rules on held-out random problems after six accepted
steps, with GD used for the first step in all methods.
From left to right the families are pseudo-Huber, subquadratic ($p=1.5$), and log-cosh. The signed recovery $R_M$, defined in
Eq.~\eqref{eq:saferecovery}, is evaluated separately
at each prescribed initial curvature ratio $\Phi_0$.
Positive (negative) values indicate a net improvement
(deterioration) relative to GD-memory.
The mean and regret maps are fitted using only the
other two families.
The binary and five-level local selectors maximize
the immediate quadratic decrease before memory
correction, while the two-step five-level oracle
provides a diagnostic reference for the potential
value of adaptation.
The vertical range is truncated at $R_M=-0.5$ to make the shallow-region behavior visible. Table~\ref{tab:policy} pools the signed gains and positive oracle gains over all 23 sampled values of $\Phi_0$, including points outside the displayed vertical range.}
 \label{fig:losorecovery}
\end{figure*}

Despite the good fit to the mean preferred exponent, the mean map gives a negative aggregate recovery in all three families, as shown in
Table~\ref{tab:policy} and Fig.~\ref{fig:losorecovery}. The
values are approximately $-0.137$ for pseudo-Huber, $-0.334$
for the subquadratic family, and $-0.128$ for log-cosh. 
These tests use 1000 new random realizations per starting-curvature value and family, with the fit for each family determined only
from the other two.

Why does the mean map perform so poorly? Consider two states with similar mismatch, one where SIGN is much better than GD and another where GD is only marginally better than SIGN. Both choices enter the average preferred exponent with the same weight. The cost of choosing the wrong exponent, however, can be quite different in the two states.

We therefore test a second rule, the \emph{regret map}. It uses the same functional form as Eq.~\eqref{eq:hmap}, but its parameters
are selected to reduce the average two-step cost penalty at the frozen training states. For an exponent in the five-level set,
we define this penalty as
\begin{equation}
 r_t(h)=
 \frac{V_t(h)-\min_{a\in\mathcal E}V_t(a)}{F(x_t)},
 \qquad h\in\mathcal E.
 \label{eq:regret}
\end{equation}
Here $V_t(h)$, defined in Eq.~\eqref{eq:Vh}, is the lowest cost
after two trial steps when the first exponent is fixed at $h$.
The penalty vanishes for a best available first choice and is
positive for more costly choices. We average these penalties in the same 24 training bins used for the mean map and interpolate linearly between the five stored exponents. The fit minimizes the population-weighted average of these interpolated regrets, evaluating the map at each bin median. 

The fitted parameters of both maps are collected
in Table~\ref{tab:mapparams}.
\begin{table}[t]
\centering
\caption{Fitted parameters of the mean and regret
maps defined by Eq.~\eqref{eq:hmap}.
For each held-out family, the parameters are
determined using only the other two families.
PH denotes pseudo-Huber, SQ denotes subquadratic
($p=1.5$), and LC denotes log-cosh.}
\label{tab:mapparams}
\small
\setlength{\tabcolsep}{3pt}
\begin{tabular}{llrrrr}
\toprule
Map & Held out & $h_{\mathrm{low}}$
 & $h_{\mathrm{high}}$ & $\log_{10}N_0$ & $\gamma$ \\
\midrule
Mean   & PH & 0.43 & 0.92 & -1.60 & 0.25 \\
       & SQ & 0.44 & 0.94 & -2.00 & 0.25 \\
       & LC & 0.43 & 0.92 & -1.60 & 0.25 \\
\midrule
Regret & PH & 0.76 & 0.99 & -0.50 & 1.60 \\
       & SQ & 0.76 & 0.99 & -0.50 & 1.80 \\
       & LC & 0.76 & 0.99 & -0.50 & 1.60 \\
\bottomrule
\end{tabular}
\end{table}
In Fig.~\ref{fig:losostruct} the regret map gives much larger exponents than the mean map. This is hardly surprising: the mean map fits the preferred exponents, while the regret map takes into account the cost of choosing a particular exponent at each state.

\begin{table*}[t]
\centering
\caption{Signed aggregate recovery $R_M$ for the
closed-loop methods on the three held-out cost-function
families. Values in brackets are $95\%$ percentile
bootstrap confidence intervals, obtained from 5000
resamples of the 1000 independent test realizations
per family. All 23 starting curvature values and
all methods are retained together within each
resampled realization. The intervals are conditional
on the fitted maps.}
\label{tab:policy}
\begin{tabular}{lccc}
\toprule
Method & Pseudo-Huber & Subquadratic ($p=1.5$)
       & Log-cosh \\
\midrule
Binary local selector
 & $-0.023$ & $-0.112$ & $-0.051$ \\
 & $[-0.093,0.039]$
 & $[-0.202,-0.028]$
 & $[-0.123,0.015]$ \\[3pt]
Five-level local selector
 & $+0.075$ & $+0.043$ & $+0.049$ \\
 & $[0.009,0.135]$
 & $[-0.041,0.115]$
 & $[-0.022,0.114]$ \\[3pt]
Mean map
 & $-0.137$ & $-0.334$ & $-0.128$ \\
 & $[-0.196,-0.083]$
 & $[-0.431,-0.249]$
 & $[-0.190,-0.071]$ \\[3pt]
Regret map
 & $-0.026$ & $-0.035$ & $-0.043$ \\
 & $[-0.054,0.001]$
 & $[-0.070,-0.003]$
 & $[-0.073,-0.015]$ \\[3pt]
Two-step oracle
 & $+0.892$ & $+0.918$ & $+0.906$ \\
 & $[0.872,0.909]$
 & $[0.898,0.934]$
 & $[0.888,0.922]$ \\
\bottomrule
\end{tabular}
\end{table*}

In complete minimization runs, the regret map also
produces much smaller losses relative to GD-memory.
The bootstrap confidence intervals in
Table~\ref{tab:policy} confirm the negative aggregate
recovery of the mean map in all three families.
For the regret map, the central estimates remain
negative in all three families, with confidence
intervals excluding zero for the subquadratic and
log-cosh families. The regret map therefore does
less harm, but neither policy reliably improves
upon GD-memory.

The two idealized targets are
\begin{equation}
   \begin{aligned}
 h_{\mathrm{mean}}(N)&=\mathbb E\!\left[\bar h_t^{(2)}\mid N_{\mathrm{sec}}=N\right],\\
 h_{\mathrm{regret}}(N)&\in\argmin_{h\in\mathcal E}
 \mathbb E\!\left[r_t(h)\mid N_{\mathrm{sec}}=N\right].
 \end{aligned}
 \label{eq:argminmean}
\end{equation}
The mean preferred exponent need not minimize the expected normalized regret. The actual fits use continuous maps and binned data. Even minimizing the expected regret at frozen states would not guarantee a better trajectory when the rule generates its own states.

For comparison we test two locally greedy selectors, one restricted to GD and SIGN and the other using all five fresh exponents. Both maximize the immediate decrease predicted by the quadratic approximation, before applying the memory correction. Ties in these local comparisons are resolved in favor of the smaller exponent.

In the test summarized in Table~\ref{tab:policy},
the five-level local selector gives a small positive
aggregate recovery in all three families.
However, the $95\%$ bootstrap confidence intervals
exclude zero only for the pseudo-Huber family.
Furthermore, the smaller independent pilot study
described in Sec.~\ref{sec:gated} yields negative
recoveries for pseudo-Huber and log-cosh.
The evidence for a systematic improvement from
local five-level selection is therefore limited.
In contrast, the large positive recovery of the
two-step oracle is well resolved in all three
families.
The five-level local selector works better in the shallow regions and deteriorates near the quadratic core.
This locally greedy selector also requires curvature evaluations for several
candidate directions, making it more expensive than either fitted map. 


\subsection{Can a curvature-based gate make adaptation safer?}
\label{sec:gated}

The performance of the five-level local selector
depends on the initial curvature regime, with a
deterioration near the
quadratic minimum. We might avoid some of the losses by using the five-level selector only where the secant--Hessian mismatch is large. We test a \emph{gated policy} that uses GD-memory unless the mismatch reaches or exceeds a chosen threshold.

The threshold determines when to use local selection; the quadratic model then determines the exponent. For a threshold $T$, we set
\begin{equation}
 h_t=\begin{cases}
  1, & N_{\mathrm{sec}}(t)<T,\\[2pt]
  \displaystyle\arg\max_{h\in\mathcal E} D\!\left(q_t(h)\right),
     & N_{\mathrm{sec}}(t)\geq T,
 \end{cases}
 \label{eq:gateclosed}
\end{equation}
where $\mathcal E=\{0,0.25,0.5,0.75,1\}$ and $D$ is the predicted
quadratic decrease defined in Sec.~\ref{sec:local}. The comparison
in the second line concerns the fresh directions, before the
memory correction. Whichever exponent is selected, the correction
and line minimization are then carried out as before. All methods
use GD for their first step, when the mismatch has not yet been
measured. Setting $T=\infty$ keeps the gate closed and reproduces
GD-memory throughout; a threshold below all observed mismatches
recovers the ungated local selector.

We determine the threshold from the outcome of complete optimization runs, rather than from the frozen-state preferences
studied in Sec.~\ref{sec:powergrid}. For each held-out family, we
try a range of thresholds on the other two families and choose the
one giving the largest mean signed gain over GD-memory, using Eq.~\eqref{eq:safegain}. For $N_{\mathrm{sec}}$ we scan the training-state quantiles $0.025,0.05,\ldots,0.95$, as well as the always-open and always-closed gates. The selected
threshold is then fixed and tested on independent realizations of
the omitted family. This leave-one-family-out test uses 60 random
problems per prescribed starting-curvature value in each training
family and 300 per value in the held-out family, with six accepted
steps per run. The search includes the possibility of never opening
the gate.

For comparison, we also construct an idealized gate using the curvature ratio $\Phi(x_t)$ instead of $N_{\mathrm{sec}}$. This gate
opens when $\Phi(x_t)$ is at or below a fitted threshold. We scan thresholds from zero to one in increments of $0.025$, together with always-closed and always-open gates. Since
$\Phi$ requires the Hessian at the known minimum, the procedure is
not generally available in practice, so it should be considered an offline control
for testing whether more direct knowledge of the curvature regime makes this switching strategy more effective. Its threshold is
trained and tested with the same separation of families.

\begin{figure*}[t]
 \centering
 \includegraphics[page=1,width=0.99\textwidth]{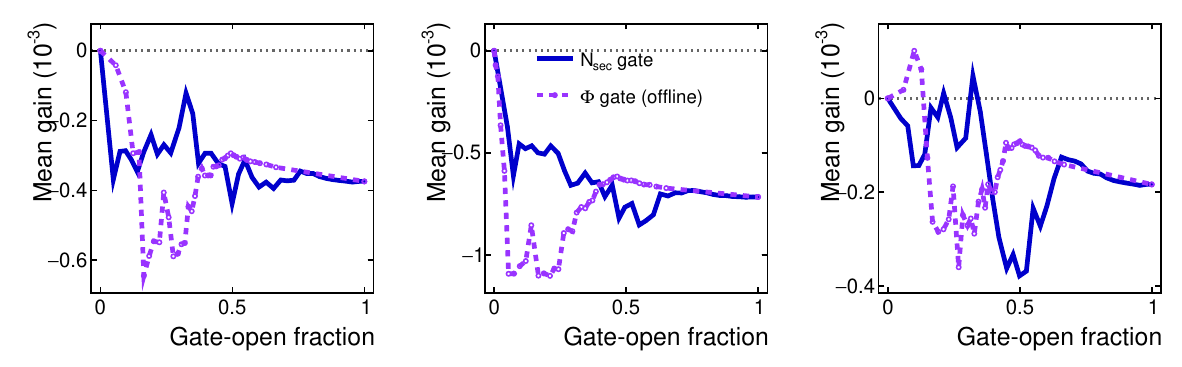}
\caption{Closed-loop tuning of a gate around the
five-level local selector. From left to right, the
held-out families are pseudo-Huber, subquadratic
($p=1.5$), and log-cosh.
Each threshold is evaluated on the other two families.
The horizontal axis gives the fraction of decisions
after the initial GD step at which the gate opens
in the training runs; the vertical axis gives the
mean signed gain $(F_6^{\mathrm{GD}}-F_6^{\mathrm{gate}})/F_0$,
in units of $10^{-3}$.
Solid curves use the secant--Hessian mismatch,
whereas dashed curves use the curvature ratio $\Phi$,
available only as an offline control.
A gate-open fraction of zero corresponds to pure
GD-memory.
In two of the three training comparisons, both gates
are optimized by never activating the local selector.
The selected thresholds are then tested unchanged
on independent realizations of the held-out family.}
 \label{fig:gatedgrid}
\end{figure*}

The threshold search is shown in Fig.~\ref{fig:gatedgrid}.
When pseudo-Huber is held out, the best threshold learned from the
subquadratic and log-cosh families is $T=\infty$: the training
procedure prefers never to activate the local selector. The same
choice is made when the subquadratic family is held out. In both
cases, the offline $\Phi$ gate also remains closed. These fitted
policies reproduce GD-memory on their respective held-out tests.
In fact, in this pilot,
the ungated five-level selector yields
$R_M\simeq-0.031$ (pseudo-Huber), $+0.199$
(subquadratic), and $-0.067$ (log-cosh) on the
held-out problems. Its positive aggregate gain
on the subquadratic family is suppressed by the
threshold chosen from the other two families,
whereas for pseudo-Huber and log-cosh the ungated
selector produces a net loss.
For these latter two families, the signs differ
from those reported in Table~\ref{tab:policy},
indicating that the positive aggregate gains
observed there are not robust across the two studies.

When log-cosh is held out, instead, training selects a finite mismatch
threshold, $T\simeq0.189$. The gate opens on approximately $35\%$
of the decisions after the initial step in the held-out tests,
reducing the number of extra curvature probes relative to the
ungated selector. Computing $N_{\mathrm{sec}}$ reuses the
Hessian--vector product needed for the memory correction; the
five additional products required to compare the local candidate
directions are evaluated only when the gate opens. The signed aggregate recovery of the mismatch gate is nevertheless $R_M\simeq-0.043$.
The offline $\Phi$ gate also selects a nontrivial threshold, but
its held-out recovery is still slightly negative,
$R_M\simeq-0.013$. Thus neither version produces a net improvement over GD-memory in this transfer test.

These are small transfer tests, especially on the training side, so we would not attach much significance to the particular threshold values. The gate trained on pseudo-Huber and log-cosh misses the positive gain of the ungated selector on the subquadratic family, while the finite threshold for log-cosh still permits unfavorable choices. Using $\Phi$ instead of $N_{\mathrm{sec}}$ does not solve the problem. For these tests we have found no way of recovering the gains of the two-step oracle with a threshold rule.



\section{Relation to previous optimization methods}
\label{sec:related}

The choice of a gradient-based search direction has long been
understood to depend on how changes in different coordinates
are measured. Ordinary gradient descent is the direction of
steepest decrease for a fixed Euclidean step length, while
sign descent is the corresponding direction when the largest
allowed coordinate displacement is fixed. The geometry of
sign descent and its relation to coordinatewise scaling have
been examined by Balles, Pedregosa, and Le Roux
\cite{Balles2020Geometry}. Sign-based updates have also been
developed for stochastic optimization, notably in signSGD,
where transmitting gradient signs reduces communication
requirements in distributed learning
\cite{Bernstein2018SignSGD}. Our experiments use deterministic
gradients and line minimizations, so communication efficiency
is not the motivation for considering SIGN. We use it as one
limiting way of weighting gradient components.

Powerball methods use a componentwise
transformation proportional to
$\operatorname{sign}(g_i)|g_i|^h$, with the exponent controlling
the influence of gradient magnitude
\cite{Yuan2019Powerball}. Related powered-gradient schemes
have been studied in stochastic optimization
\cite{Zhou2020PbSGD}. These methods demonstrate that
compressing gradient magnitudes can improve optimization
under suitable conditions. In our study, the same family
of transformations provides a controlled set of possible
actions. The exponent is allowed to change from one
iteration to the next, and its consequences are evaluated
together with the curvature-dependent memory of the
preceding direction.

A different and widely used form of adaptation changes the
effective step size of individual coordinates using the
history of their gradients. AdaGrad and Adam are familiar
examples \cite{Duchi2011AdaGrad,Kingma2015Adam}. Such
methods can compensate for differences in parameter
sensitivity, but their mechanisms differ from changing a
single exponent applied to the current gradient. In our
construction, the coordinates are fixed and scaled in
advance, while the exponent controls the relative weight
of the present gradient components. The resulting direction
is subsequently modified using curvature information from
the preceding step. These choices allow us to isolate how
changes in gradient geometry interact with directional
memory.

A similar issue arises in low-rank gradient optimization.
GaLore periodically refreshes the subspace in which
gradient statistics are stored, so changing the projection
can make information accumulated in the previous
coordinates less useful~\citep{Zhao2024GaLore}.
LDAdam and SubTrack++ address this problem by
adapting the stored optimizer statistics when the
subspace changes~\citep{Robert2025LDAdam,Rajabi2025SubTrack}.
These methods work in a rather different setting,
but they face a related problem: changing the way
an update is constructed also affects how information
from earlier steps should be used.

The memory correction is closely related to conjugate-gradient
methods. For a positive-definite quadratic cost function,
conjugate directions avoid repeatedly undoing progress
along previously explored directions, and classical
conjugate gradients exploit this structure very effectively
\cite{HestenesStiefel1952,NocedalWright2006}. 
Our correction enforces conjugacy with the most recent
search direction using a Hessian--vector product, and
can be applied after any of the chosen gradient
transformations. As discussed in Sec.~III, retaining
GD as the fresh direction throughout the trajectory
recovers linear conjugate gradients on a positive-definite
quadratic under exact line minimizations.
Changing the fresh gradient transformation removes
the guarantee of global conjugacy, although conjugacy
with the preceding direction is still enforced.
A related phenomenon arises in conjugate gradients
with variable preconditioning, where the advantage
over steepest descent can disappear even when the
preconditioners are positive definite and the
condition numbers of the preconditioned systems
are uniformly bounded~\citep{knyazev2007variable}.
Changing the fresh gradient transformation plays
an analogous role in our construction.

Differences between gradients at consecutive points are central to
secant conditions in quasi-Newton optimization
\cite{NocedalWright2006}. Our secant--Hessian mismatch
compares the observed change in gradient with the change
predicted by applying the Hessian at the endpoint of the
step. We use the discrepancy to detect curvature changes along the step. The required Hessian--vector product is already used by our memory correction.
A related diagnostic of the validity of a local
quadratic model is the trust-region ratio of actual
to predicted cost decrease~\citep{NocedalWright2006}.
The trust-region ratio compares cost reductions using a model at the beginning of the step; our diagnostic compares gradient changes using the Hessian at the end. It would be interesting to see which is a better predictor of the preferred exponent.

Learning-to-optimize methods seek to infer update rules from experience on related minimization problems
\cite{Andrychowicz2016Learning}. Our setting is more
restricted: the candidate directions are prescribed, and
the two-step searches serve as numerical references for
judging the value of individual choices. This restriction
allows us to compare a favorable action at a given state
with the behavior of a rule that repeatedly selects
actions during minimization.
A related distinction between short- and long-term
performance has been identified in stochastic
meta-optimization, where tuning optimizer
hyperparameters over short horizons can lead to
unfavorable choices over longer
runs~\citep{Wu2018ShortHorizon}.
From a control-theoretic perspective, our two-step
oracle implements a form of model-predictive control
with an exact model of the optimizer's state
transitions~\citep{RawlingsMayne2017}.
Even in the absence of model uncertainty, optimizing
over a restricted prediction horizon does not
guarantee optimality over a longer trajectory.

Here we can compare the immediate gain from a direction with the gain over two or more steps while retaining the same memory correction. The quadratic example isolates the effect in an exactly solvable setting, and the nonlinear tests show how difficult it can be to turn the potential gains into a useful rule.

\section{Discussion}
\label{sec:discussion}

Testing an adaptive rule on fixed states is useful, but one should also follow the trajectories generated by the rule itself. This is especially important with directional memory: every decision changes both the point where the next gradient is evaluated and the direction against which the next step is corrected. In our experiments a good description of the average preferred exponent did not suffice to improve the complete minimization.

A single number may be too little information to decide how to change the gradient geometry.
The secant--Hessian mismatch tells us how well the Hessian
at the end of a step accounts for the gradient change along
that step, but it does not specify how the curvature is
distributed among different directions, nor how the
available choices will affect the next remembered direction.
Consequently, states with similar mismatch can have quite
different costs associated with choosing a particular
exponent. We might do better by estimating the possible gains and losses from departing from GD-memory than by trying to reproduce the average exponent selected by the oracle.

We would also like to understand more about the role of memory.
Our exact quadratic example establishes that an immediate
improvement can produce an unfavorable two-step outcome
under the specified correction rule. In the nonlinear
experiments, however, a change of direction also changes
the subsequent position, so the contributions of the new
position and the remembered direction cannot be separated
directly from the performance curves. Controlled comparisons
with alternative memory updates, or with memory disabled,
would help determine which features of the observed behavior
arise specifically from the curvature-conjugacy constraint.
Such experiments would also clarify the extent to which the
deterioration near quadratic minima is associated with
disrupting a favorable sequence of directions.

There is also the question of computational cost. The secant--Hessian mismatch can be evaluated
with curvature information already used by our memory
correction, whereas testing several candidate directions
requires additional Hessian--vector products. The two-step
oracle is still more demanding, since it evaluates
speculative trajectories that are never executed. We have compared methods at a fixed number of accepted steps; an efficiency comparison would need to account for the additional evaluations. A practical implementation
would need to obtain improvements large enough to justify
the information required to select each step.

Our benchmarks are quite restricted. The three cost-function families
share a common construction, have positive curvature, and
are studied in six dimensions using accurate gradients,
Hessian--vector products, and line minimizations.
The step budget of six also coincides with the
finite-termination bound for linear conjugate gradients
in this dimension. The relative performance of the
methods in the near-quadratic regime may therefore
depend on the chosen step budget and dimensionality.
Their
known minima allow us to classify the landscapes by a
curvature ratio that would generally be unavailable in an
actual fitting problem. Moreover, both power-gradient
directions and the secant--Hessian mismatch depend on the
scaling of the coordinates. We have kept that scaling fixed
to permit meaningful comparisons, but a method intended
for physical applications would need to be examined under
different parameterizations, in higher dimensions, and in
the presence of numerical or measurement noise. The
agreement found among our three families should therefore
not be interpreted as establishing a universal relation
between curvature variation and the preferred gradient
exponent.

The two-step search provides examples where departures from GD-memory help and examples where they hurt. We would like to find a small set of readily available quantities that distinguish these cases, taking into account the size of the possible gain or loss. Any new rule will have to be tested in complete minimizations, including the computational cost of choosing each step.

\section{Conclusion}
\label{sec:conclusion}

We have studied the effect of changing the gradient geometry in a minimization algorithm with curvature memory. Even for a positive-definite quadratic, choosing the direction giving the largest immediate decrease can lead to a higher cost after two steps. On three nonlinear benchmarks we found that changing the exponent can be useful in the shallow regions, and that intermediate powers are often preferred in the two-step comparisons.

The secant--Hessian mismatch is easy to evaluate
in our algorithm and gives a good prediction of
the mean preferred exponent across the three
families. Yet neither of the fitted maps nor the
tested threshold rules reliably improves on
GD-memory. The lesson from these results is that
knowing which direction works best on average
does not tell us which one to choose at a particular
step. With curvature memory, that choice also
affects the directions available at later steps.

\appendix

\section{Additional details of the three-dimensional counterexample}
\label{app:counter}

For the example in Sec.~\ref{sec:greed}, the remembered
direction is $r=(0,16,-1)$. With the Hessian of
Eq.~\eqref{eq:counterH}, we have
\[
 Hr=(0,64,-16),\qquad r^THr=1040.
\]
The memory correction of Eq.~\eqref{eq:memory}
therefore gives $\beta_G=-12/65$ for GD and
$\beta_S=3/65$ for SIGN. The corresponding
corrected directions are
\begin{align}
 p_G&=\left(-a,-\frac{257}{65},-\frac{1028}{65}\right),\\
 p_S&=\left(-1,-\frac{17}{65},-\frac{68}{65}\right).
\end{align}
The fresh and projected decreases are
\begin{align}
 D_G^{\rm fresh}(a)&=\frac{(a^2+257)^2}{2(a^2+4100)},\\
 D_S^{\rm fresh}(a)&=\frac{(a+17)^2}{42},\\
 D_G^{\rm proj}(a)&=\frac{65(a^2+257)^2}{2(65a^2+264196)},\\
 D_S^{\rm proj}(a)&=\frac{65(a+17)^2}{2442}.
\end{align}
For the relevant two-step branches,
\begin{align}
 f_{GG}(a)&=
\frac{66770568\,a^4}{4330625a^4+2358873986a^2+283560576065},\\
 f_{SG}(a)&=\frac{2(68a-65)^2}{557801}.
\end{align}

Over the range $0.1\leq a\leq3$, we have verified
numerically that GD gives the better second step
for either first-step choice, unless the minimum
has already been reached.
The two-step comparison therefore reduces to the
expressions for $f_{GG}(a)$ and $f_{SG}(a)$ above.
The boundaries quoted in Sec.~\ref{sec:greed}
are obtained by equating the corresponding
analytical expressions:
\begin{align*}
D_S^{\rm proj}(a)=D_G^{\rm proj}(a)
    &:\quad a\simeq0.4870,\\
f_{SG}(a)=f_{GG}(a)
    &:\quad a\simeq0.8667,\;1.0991,\\
D_S^{\rm fresh}(a)=D_G^{\rm fresh}(a)
    &:\quad a\simeq1.5620.
\end{align*}

\section{Supporting experimental controls}
\label{app:controls}

We briefly report some of the additional tests that helped guide the main study.

\paragraph{Binary versus continuous geometry with memory.}
In preliminary experiments on rotated nonlinear cost
functions, retaining one curvature-conjugate direction
often produced larger improvements than refining the
choice of gradient exponent. These observations motivated
the controlled comparison of binary and intermediate
gradient transformations in the main text.

\paragraph{Memory-aware greedy $h$.}
Additional exploratory tests examined the selection
of the exponent after the memory correction. Such
selection could produce worse trajectories, and
increasing the resolution of the local search did not
necessarily improve the final result. These observations
motivated the distinction between selecting a fresh
direction and evaluating its consequences after the
memory correction, discussed in
Secs.~\ref{sec:memory} and \ref{sec:greed}.


\paragraph{One-HVP local diagnostics.}
Several recycled quantities were screened as
predictors of the preferred first action in the
binary two-step search, including the gradient
evenness factor $G$ defined in Sec.~\ref{sec:local}
and the relative squared size of the memory correction, $\|\beta_t p_{t-1}\|^2/\|q_t\|^2$.
None of these quantities gave enough information to recover the gains of the oracle. The secant mismatch $N_{\mathrm{sec}}$ was the most robust indicator of departures from quadratic behavior; it was much less useful for choosing an individual step.

\bibliographystyle{apsrev4-2}
\bibliography{references}

\end{document}